
\input amstex

\documentstyle{amsppt}

 \magnification\magstep1
\NoBlackBoxes
\NoRunningHeads

\voffset=0.8cm
 \hsize30pc
 \vsize46pc

\define\Z{\Bbb Z}

\define\Q{\Bbb Q}
\define\C{\Bbb C}
\define\F{\Bbb F}

\define\calC{\Cal C}

\define\Aut{\text{\rm Aut}}

\define\k{\ovl{k}}

\define\Imag {\text{\rm Im}}

\define\lra{\longrightarrow}

\define\ovl{\overline}
\define\BP{\Bbb P}

\topmatter
\title  d-gonality of modular curves and   \\ bounding torsions
\endtitle

\author  Khac Viet Nguyen
\& Masa-Hiko Saito
\endauthor
\thanks The first author was supported by JSPS through the long
term visiting program (L-95504); \endthanks

\abstract
We study the problem of $d$-gonality of the modular curve $X_0(N)$.
 As a result,
we can give an upperbound of the level $N$ by means of $d$.
 This generalizes
Ogg's result on hyperelliptic modular curves ($d = 2$) in ([O]).
As a corollary of this result,
we  prove an analogue of the strong Uniform Boundedness Conjecture
for elliptic curves defined over the function fields of curves.
If a base curve is $d$-gonal, we can bound orders of torsions of
Mordell-Weil groups in terms of $d$ uniformly.
\endabstract

\affil Hanoi Institute of Mathematics \&
 Kyoto University
\endaffil
\address  Hanoi Institute of Mathematics, 631, Bo Ho, 10000,
Hanoi, Vietnam  \quad \&
\linebreak
 Department of Mathematics, Faculty of
Science, Kyoto  University, Kyoto,  606-01, Japan
\endaddress

\email nkviet\@thevinh.ac.vn  \quad \&  \quad
 viet\@kusm.kyoto-u.ac.jp
\endemail

\address Department of Mathematics, Faculty of Science, Kyoto
University, Kyoto,  606-01, Japan
\endaddress

\email mhsaito\@kusm.kyoto-u.ac.jp
\endemail

\endtopmatter

\noindent{\bf \S 0 Introduction}

\vskip 0.5cm

Let $N$ be a positive integer, and let
$\Gamma_0(N)$ be the subgroup of $\Gamma= SL_2(\Z)/ (\pm 1)$
defined by matrices $\left(  \matrix a & b \\ c & d \endmatrix
\right)$ with $N$ dividing $c$.    Then $\Gamma_0 (N)$ acts on
the upper half plane
${\Bbb H}$ properly  discontinuously, and let
$$
Y_0(N)_{\C} = \Gamma_0(N) \backslash {\Bbb H}.
$$
The modular curve  $X_0(N)_{\C}$ is the compactification of
$Y_0(N)_{\C}$ obtained by adding  the cusps. A smooth projective
curve $C$ defined over an algebraically closed  field $k$ is
called {\it d-gonal} if there exists a  finite morphism $f:X
\lra \BP^1_k$ of degree $d$.    For example, if a smooth curve
$C$ is 1-gonal, $C$ is isomorphic to
$\BP^1_k$, and if $C$ is 2-gonal, then either $g(C) \leq 1$ or
$C$ is a hyperelliptic curve of genus $g(C) \geq 2$.

\vskip 0.2cm
In this paper, we shall consider the following
problem.

\proclaim{Problem 0.1}
 If $X_{0}(N)_{\C}$ is $d$-gonal, can one
give a bound for $N$  by means of $d$?
\endproclaim

\vskip 0.2cm By the genus formula of $X_0(N)_{\C}$, one can
determine  the case when $g(N) :=g(X_0(N)_{\C}) =  0$, which
gives the answer  for $d = 1$.  (One knows that $ N \leq 25$ for
$d = 1$.) For $d = 2$,  Ogg [O] classified all  hyperelliptic
modular curves $X_{0}(N)$,  and it gives us the sharp  bound $N
\leq 71$.   In this paper, we deal with  Problem  0.1 for
general $d$.

Our main theorem can be stated as follows.

\proclaim{Theorem 0.2}
(cf. Corollary 3.5.)
Let $X_0(N)$ be the
modular curve of level $N$.  If $X_0(N)_{\C}$  is d-gonal,
  then we have
$$
 N  \leq   \left\{ \matrix  25   & \text{if $d = 1$}
 \\   71   & \text{if $d = 2$}  \\
(48(d - 1)^2 + 35)( 36(d -1)^2 + 47 )
& \text{if $d  \geq 3$}
\endmatrix \right. . \tag 0.1
$$
\endproclaim

Note that we can obtain better bound for $N$ odd (see Corollary
3.3).

\vskip 0.2cm

Before going into the sketch of the proof of Theorem 0.2,   we
explain the motivation of  Problem 0.1.   Our original
motivation is a function field analogue of the so-called  strong
Uniform Boundedness Conjecture (abbreviated as sUBC). According
to [Kam-Ma],  we denote by $\Phi(d)$  the set of all isomorphism
classes of  finite abelian groups occurring as the full groups of
torsion  in the Mordell-Weil groups of elliptic curves over
number  fields  $K$ with absolute degree $[K:\Q] \leq d$. Then
sUBC states that $\Phi(d)$ is finite.

Let $d$ be a positive integer.  A prime number $p$ is called a
{\it torsion prime for degree $d$} if there is a number field
$K$ of  degree $d$, an elliptic curve $E$ defined over $K$, and a
$K$-rational  point $P$ of $E$, of order $p$.  Let $S(d)$ denote
the set of  torsion primes of degree $\leq d$.  Then Kamienny and
Mazur [Kam-Ma]  showed that $\Phi(d)$ is finite if and only if
$S(d)$ is finite  (see also [E]). Recently, Merel [Me] showed
that for $p \in S(d)$
$$
 p < d^{3d^2},
$$  which completes the proof of sUBC.  (Note that one could not
give an explicit bound for $N$ occurring as an exponent of a
 group in $\Phi(d)$.  The complete descriptions of $\Phi(d)$
for   $d = 1, 2$ is known ( see  [Ma], [Kam-Ma],  [E] and references
therein).

\vskip 0.2cm In this paper, we prove  the following theorem which
can  be viewed as an analogue of sUBC in the function field case.
Let $C$ be a smooth projective curve defined over the field of
complex numbers and set $K = \C(C)$,  the function field of
$C$.
Let  $E$ be a non-constant elliptic curve defined over $K$.
 If $C$ is $d$-gonal, then
 there exists an extension of  fields $\C(\BP^1_{\C}) = \C(t)
\hookrightarrow K$  of degree $d$.  (Note that such a  field
extension may not
 be unique even if  we assume that $d$ is minimal.)

Now we denote   by $\Phi_{fun}(d)$
 the set of all isomorphism classes of  finite abelian groups
occurring as the full groups of torsion  in the Mordell-Weil
groups of non-constant   elliptic curves over a function field
$K$ with an extension $\C(\BP^1) \hookrightarrow K$ of degree
$\leq d$.  Then  we propose an function field analogue of sUBC
as:  {\it $ \Phi_{fun}(d)$ is finite.}

This seems to be a quite natural question, which originally
arose in a seminar at Kochi (July, 1995) with H. Tokunaga and
the second author.   Furthermore  they realized that the analogue
can be reduced to the problem of $d$-gonality of modular
curves.

\vskip 0.2cm
{}From Theorem 0.2, we obtain the  following theorem,
which obviously implies that
$ \Phi_{fun}(d)$ is finite.

\proclaim{Theorem 0.3 (An analogue of sUBC for function field
case)} (cf. Theorem 4.3.)
 Let $E$ be a non-constant elliptic curve defined over  the
function field $K$ of a complex smooth projective curve $C$.
Assume that $C$ is d-gonal and the Mordell-Weil group
$E(K)$ has a torsion element of order $N$.  Then
$$ N \leq B(d)
$$ where $B(d)$ is the right hand side of the inequality  (0.1).
\endproclaim

Let us give a short summary of this paper.    In \S 1, we recall
necessary results on curves.  One important  remark is that
$d$-gonality  of a curve descends to a curve  below via a finite
morphism, and another is the so-called
 inequality of Castelnuovo -Severi.
In \S 2, we will prove
``Tower Theorem 2.1'', which is one of  our main contributions in
this paper.  In \S 3, we will  generalize a beautiful argument of
Ogg [O] and Harris-Silverman [Ha-Si]  to obtain a bound of
level $N$ by means of $d$.
{}From Tower Theorem 2.1, we obtain a
finite morphism $f:X_{0}(N)_{\Q} \lra C'$ defined over $\Q$ of
degree $d' \leq d$  such that $g(C') \leq (d/d' -1)^2$.  For a
prime $p$ not dividing
$N$,
$X_{0}(N)_{\Q}$ has a good reduction at $p$.
When $g(C')  \geq 1$,  applying
Good Reduction Lemma 5.1, $C'$ has also good reduction
at $p$ and we obtain a  finite morphism
$f_{\F_{p}}:X_{0}(N)_{\F_p} \lra  C'_{\F_p}$ of degree $d$
defined over $\F_p$.

 By Weil's theorem of the analogue of the Riemann
hypothesis, we can bound the number of $\F_{p^2}$-rational
points of $ C'_{\F_p^2}$ from the above by $g(C')$, hence
by $(d/d' -1)^2$.  Therefore the number of
$\F_{p^2}$-rational  points of $X_{0}(N)_{\F_p^2}$
is bounded from above only by $d$.  On the other hand,
we have enough $\F_{p^2}$-rational
points of $X_{0}(N)_{\F_p^2}$ from  points
corresponding to supersingular elliptic curves  and
cusps.  For example, if $p =2$, $\#(X_{0}(N)_{\F_4})$
is at least $\frac{1}{12} (N + 1) + 2$.
This argument gives us bound of $N$ by means of $d$ and
a proof of Theorem 0.2.
  In \S 4, we prove Theorem 0.3 by using Theorem 0.2 and  the moduli
property of $X_{0}(N)$.  In \S 5, we shall prove  the Good
Reduction Lemma 5.1 which is  our another contribution in this
paper.

\vskip 0.5cm

Let us mention some  results  related  to ours.   In the case
of function field $K = k(C)$ with  char $k = 0$,  a bound of
$|E(K)_{tors}|$ in terms of the genus $g = g(C)$   was  given
first by Levin [Le],  and later (with shaper  bound) by Hindry
and Silverman [Hi-Si].   A rough estimate  in    [Hi-Si, Theorem
7.2, (c)] gives us
$$
|E(K)_{tors}| \leq 144 (g + 1)^{2/3}.
$$
(See [Hi-Si, Theorem 7.2, (a),(b)] for refined bounds.)
The bound in [Hi-Si, Theorem 7.2, (a)] follows from  the function
field analogue of Szpiro's conjecture ([Hi-Si, Theorem 5.1]),
which was originally   proved by Kodaira and Shioda (cf. [Shi,
Proposition 2.8]).   Note
that since the genus $g$ can be arbitrarily large with
$d$ fixed, so it is still far from  the uniform bound by $d$.

As far as  the set $\Phi_{fun}(d)$ is concerned,  Cox-Parry [C-P]
determined the set $\Phi_{fun}(1)$ completely. (See also [M-P]).
We can also determine the set $\Phi_{fun}(2)$  completely.
(See Nguyen's report of our joint work [Ng]).

When  $k$ is a finite field with $q$ elements of  char $k = p >
0$, we remark that
 our method in  this paper gives a bound of order of {\it
prime-to-$p$-part of}
 $E(K)_{tors}$  by means of $d$.   Our  bound inevitably depends
on $d$ and the characteristic
$p$ but not depend on $q$.   Recently,  Goldfeld and Szpiro gave
a bound for $|E(K)_{tors}|$ when
$k$ is a finite field with $q$ elements (cf. [G-S, Theorem 13]).
Their bound depends on $q$,  the genus $g$ and the inseparablilty
degree of the $j$-function $C \lra \BP^1$ associated to  a
minimal model of $E/K$.

\proclaim{Ackowledgement} {\rm We would like to thank Professor
Shigefumi Mori for  discussions  on minimal models, which have
been  very helpful for our proving  the good reduction lemma. We
are also grateful to Professor Frans  Oort who kindly answered
our question on the stable reduction  of curves via email.    The
second author would like to thank Professor  Hiro-o Tokunaga  for
discussions
 on this problem at Kochi University  which is a starting point
of this work. We are also grateful to Professor Tetsuji Shioda
for giving us kind comments and
informing  us  of some of references related to our work.

The first author would like to thank JSPS for  their financial
support which made possible this  joint work at Kyoto
University.}
\endproclaim

\vskip 0.8cm

\noindent{\bf \S 1 Some results on curves}

\vskip 0.5cm In this section, we  recall some known  results  on
curves which we  use later.   Let $k$ be an algebraically closed
field of any  characteristic $p \geq 0$.

\proclaim{Definition 1.1} {\rm A smooth projective curve $C$
defined over  $k$ is called  {\it $d$-gonal} if there exists a
finite  morphism $f:C \lra \BP^1_k$ of degree $d$, or
equivalently,  if $C$ has a $g^1_d$.   We call such  a map $f$ a
{\it $d$-gonal map}. }
\endproclaim

If $C$ admits a $d$-gonal map $f:C \lra \BP^1_k$,  then it
induces  a  field extension
$k(\BP^1) = k(t) \hookrightarrow K = k(C) $ of degree $d$, and
conversely such a field extension gives a $d$-gonal map $f:C \lra
\BP^1_k$.

\proclaim{Remark 1.2} {\rm  One may like to use the term
``$d$-gonality" for
 the case with $d$ minimal.  However the definition above  allows
us to state our results  in  simpler form and the reader can
easily apply the results to the case with $d$  minimal.
Moreover we remark that a $d$-gonal map $f:C \lra \BP^1_k$  may
not be unique up to automorphism of $C$ and $\BP^1_k$  even if
$d$ is minimal.

  Moreover,  for $d \geq (1/2) g +1$,  any curve of genus $g$
 is   $d$-gonal, and  for $d < (1/2) g  + 1$, there exist curves
 of genus $g$ having no $d$-gonal map.
  On the other hand,  there exists a  hyperelliptic curve of
genus $g$ (2-gonal curve) for any $g \geq 2$. }
\endproclaim

First of all, we recall the following lemma.  Though it  may be
well-known to experts,  for completeness,   we will give a proof
which is valid for any algebraically  closed base field $k$.
(cf. [Ne, Theorem VII.2]).

\proclaim{Lemma 1.3} Let $h:C_1 \lra C_2$ be a finite  morphism
between  two smooth projective  curves.   If
$C_1$ is d-gonal, then $C_2$ is also $d$-gonal.
\endproclaim

\demo{Proof} Let $f: C_1 \lra \BP^1$ be a rational function of
degree $d$.   Let $K_i, ( i = 1, 2)$ denote the rational function
fields of
$C_i (i =1, 2)$.  Then we have a field extension
$K_2 \hookrightarrow K_1$ of degree $m$.  Take an algebraic
closure
$\ovl{K_1}$  of $K_1$, and let $l $ and $s$ be the separable and
inseparable degrees of the extension $K_1/K_2$.   (When char of
$k$ is $0$, $s =1$ and  when char $k = p > 0 $, then $ s = p^e$
for some $e \geq 0$.) For any  element $h \in K_1$, we can
define  the norm
$N_{K_1/K_2}( h)$ as
$$ N_{K_1/K_2} (h) = \prod_{i =1}^{l} (\sigma_i(h^s))
$$ where $\sigma_i:K_1 \hookrightarrow \ovl{K_1}, {i =1,\cdots,
l}$ are different  $K_2$-embeddings  of $K_1$ into $\ovl{K_1}$.
Then   the norm $N_{K_1/K_2} (f)$ gives a  morphism
$C_2 \lra \BP^1$ of degree $\leq d$.   It is easy to see that for
some constant $\lambda \in k $  the norm
$ N_{K_1/K_2}(f - \lambda) $ gives a finite morphism
$C_2 \lra \BP^1$  exactly of degree $d$.
\enddemo

\quad
\quad  The next important result,  which we use in  the proof of
Tower Theorem 2.1,  is the so-called inequality of
Castelnuovo-Severi.

\proclaim{Lemma 1.4} ({\it The inequality of
Castelnuovo-Severi}).  Let
$C, C_1, C_2$ be  smooth projective curves over $k$ of genera $g,
g_1, g_2$ respectively
 and let   $\pi_1:C \lra C_1$ and $\pi_2:  C \lra C_2$
surjective morphisms of degrees
$d_1, d_2$ respectively. Assume  that the induced map
$$
\pi_1 \times \pi_2:C \lra C_1 \times C_2 \tag 1.1
$$  is birational onto its image.   Then one has

$$  g \leq d_1g_1+d_2g_2+(d_1-1)(d_2-1).  \tag1.2
$$
\endproclaim

\demo{Proof}  Though there are many proofs of this inequality
(cf. [H, Ch\. V],  [Gr]),  we will  give a proof using genus
formula and the Hodge index theorem.
 Set $ S = C_1 \times C_2$ and
$D = \pi_1 \times \pi_2(C) \subset S$.   Then one can calculate
the virtual genus
$p_a(D)$ of a curve $D$ on a surface $S$  as $p_a(D) = (K_S +
D)\cdot D/2 + 1$.  Then by assumption that
$C \lra D$ is birational, we have $g(C) \leq p_a(D)$.   On the
other hand, from the assumption, we obtain
$K_S \cdot D = (g_1 -1) d_1 + (g_2 -1) d_2$ and the Hodge index
theorem implies that $D^2 \leq 2 d_1 d_2$ (see [H, Ex. 1.9, Ch\.
V]).   Therefore we obtain the inequality.
\enddemo

\vskip 1cm

\noindent{\bf \S 2  Tower theorem }

\quad

In this section, we  prove the following theorem which we  call
``Tower theorem''.
Let $k$ be a  perfect field and fix  an algebraic closure $\k$
 of $k$ and let $G_{k} = Gal(\k/ k) $ be the  Galois group of the
extension $\k/k$.

\proclaim{Theorem 2.1} Let $C$ be a projective smooth curve
defined over a perfect field
$k$  and $f:C \lra \BP^1_{\k}$ a dominant morphism defined over
$\k$ of degree $d$.  Then there exists a smooth  projective curve
$C'$ defined over $k$ and a dominant morphism
$$  f':C \lra C'
$$  defined over $k$  of degree $d' $ dividing $d$ such that
$$  g(C') \leq (d/d' -1)^2. \tag 2.1
$$
\endproclaim

\demo{Proof}  At first,  we consider the following
\newline

{\bf Claim:} There exists a tower of (projective smooth) curves
over $\k$
$$  C  \lra C_n \lra \cdots \lra C_2 \lra C_1 \lra C_0 \simeq
\BP^1 \tag 2.2
$$  satisfying the following conditions:

(i) $h_i:C_i \lra C_{i-1}$ and $f_n:C \lra C_n$ are  morphisms  of
degree $e_i \geq 2$  and of degree $d_n$ defined over $\k$.

(ii) Set  $f_i:C \lra C_i$ and $d_i = \deg f_i$.  Then we have
$ d = d_0 > d_1 > \cdots > d_n \geq 1$ and for any $0 \leq i \leq
n$
$$  g(C_i) \leq (d/d_i - 1)^2.
$$

(iii) For any $\sigma \in G_{k}:= \text{\rm Gal}(\k/k)$,  the
morphism
$$  f_n \times f_n^{\sigma}:C \lra C_n \times C_n^{\sigma}
$$  has degree $d_n = \deg f_n$ onto its  image.

Let us first explain how the theorem follows from this claim.
Consider the morphism $f_n \times f_n^{\sigma}:C \lra C_n \times
C_n^{\sigma}$ for  any $\sigma \in G_k$ and let $D  = f_n \times
f_n^{\sigma}(C)$ be  the image.  Let  $p_1:D \lra C_n$ and
$p_2:D \lra C_n^{\sigma}$ be the natural projections.
  Since  $f_n$ and $f_n^{\sigma}$ factor through
$ C \lra D $ and $p_1$ and $p_2$ respectively, from  the
assumption (iii), we know that
$\deg p_i = (\deg f_n )/(\deg (C \lra D)) =  d_n/d_n = 1$,  which
implies that  the projection maps $p_1$ and
$p_2$  are birational. Hence we infer that
 $C_n \simeq C_n^{\sigma}$ (over $\k$) for  any $\sigma \in
G_k$.  This shows that $C_n$ is isomorphic (over $\k$)  to a
curve $B$ defined over $k$. Replacing $C_n$ by this new curve
$B$,  we obtain a morphism $h: C \lra B$ of degree $d_n$.  Though
$h$ may not  be defined over $k$, for any $\sigma \in G_{k}$,  the
argument given above shows that $h:C \lra B$ and
$h^{\sigma}:C \lra B$  differ by an automorphism of $B$.  That
is,  we have  $\alpha_{\sigma} \in
\Aut(B \otimes \k)$ such that
$$
 h^{\sigma} = \alpha_{\sigma} \circ h.
$$   It is easy to see that $\alpha:G_{k} \rightarrow
\Aut(B \otimes \k)$ defines a cocycle in
$H^1(G_{k}, \Aut(B \otimes \k))$.    Then by general theory of
twisting (cf. [Si., X, Theorem 2.2]), $\alpha$ corresponds  to a
twist , that is, there exists a smooth  curve
$C'$ defined over $k$ and $\k$-isomorphism
$\lambda:C' \lra B$ such that
$\alpha_{\sigma} = \lambda^{\sigma} \lambda^{-1}$.   Now the map
$\lambda^{-1} \circ h:C \lra C'$ gives a dominant morphism defined
over $k$ of degree $d_n$ and $g(C') = g(B) = g(C_n)
\leq (d/d_n -1)^n$, which implies the theorem.

Now we  show how one can obtain a tower  in  the claim.

   Take $\sigma \in G_{k}$,  and consider the map
$$
 f \times f^{\sigma}:C \lra \BP^1_{\k} \times \BP^1_{\k}.
$$  Let  $D = f \times f^{\sigma}(C) \subset
\BP^1 \times \BP^1$  be the image.
   If for any $\sigma \in G_{k}$,
$\deg  f \times f^{\sigma}  = \deg f = d$,
then  we can take
$C_n = C_0 (\simeq \BP^1_{\k})$ as the normalization of $D$ and
$f_n = f_0:C \lra C_0$ as the natural morphism to obtain  the
desired tower.
(Note that $d_n = d$ and $g(C_0) = 0$.) Hence in
this case we get a tower (2.2) as desired.

Otherwise, for some $\sigma \in G_k$, the degree of
$f \times f^{\sigma}$  onto its image is $d_1 < d = d_0$.  Then setting
$C_1 = $ the normalization of the image $f \times f^{\sigma}(C)
\subset \BP^1 \times \BP^1$, we obtain finite morphisms
$$
h_1:C_1 \lra C_0 \simeq \BP^1_{\k}, \quad h_1^{\sigma}:C_1 \lra C_0 \simeq
\BP^1_{\k},
\quad f_1: C \lra C_1
$$
such that $ f =    h_1 \circ f_1$,  $f^{\sigma} = h_1^{\sigma} \circ f_1$ with
$e_1 = \deg h_1 = \deg h_1^{\sigma} \geq 2$ and $d_1 =  \deg f_1 = d/e_1 < d$.
 Then since $h_1 \times h_1^{\sigma}: C_1 \lra \BP^1_{\k} \times \BP^1_{\k}$ is
birational onto its image,
 the inequality of Castelnuovo-Severi implies that
$$ g(C_1) \leq (e_1 - 1)^2 = (d/d_1 - 1)^2.
$$

\quad After continuing these procedures,  we may assume  that we
have a tower of (smooth)  curves
$$  C  \lra C_i \cdots  \lra C_1 \lra C_0 \simeq \BP^1
$$ satisfying only the conditions (i), (ii) up to a level $i > 0$.
For any $\sigma \in G_{k}$ and the given morphism $f_i:C \lra
C_i$,
 consider the morphism
$$  f_i \times f_i^{\sigma}:C \lra C_i \times C_i^{\sigma}.
$$

If for all $\sigma \in G_{k}$ $ f_i \times f_i^{\sigma}$ is
birational onto its image, then the tower also  satisfies the
condition (iii).  Then we can stop the procedure.

Otherwise,  for some  $\sigma \in G_k$ , the  degree of $f_i
\times f_i^{\sigma}$ is $d_{i+1} < d_i$.  Then  again let
$C_{i+1}$   be the normalization of the image of $f_i \times
f_i^{\sigma}(C)$ and  let $f_{i+1}:C \lra C_{i+1}$  be the induced
morphism  and $C_{i+1} \lra C_i \times C_i^{\sigma}$ the induced
 map birational onto its image.   Since  the degree of  each
projection $C_{i+1} \lra C_i$ and $C_{i+1} \lra C_i^{\sigma} $
is $e_{i+1} = d_i/d_{i+1} \geq  2$,  from the  inequality of
Castelnuovo-Severi  (Lemma 1.4)   and the assumption
$g(C_i) \leq (d/d_i -1)^2$, we obtain:
$$  g(C_{i+1}) \leq 2 e_{i+1} \cdot  g(C_i) + (e_{i+1} -1)^2
\leq 2 e_{i+1} \cdot (d/d_i -1)^2 + (e_{i+1} -1)^2. \tag 2.3
$$  Since $d_i = d_{i+1} \cdot e_{i+1}$, we obtain $d/d_{i+1} =
(d/d_i) \cdot e_{i+1}$, and since
 $e_{i+1} \geq 2$, we can easily see that
$$  2 e_{i+1} (d/d_i -1)^2 + (e_{i+1} -1)^2 \leq
 ((d/d_i) \cdot e_{i+1} -1)^2 =(d/d_{i+1} - 1)^2. \tag 2.4
 $$
 Hence,  together with (2.3),  we obtain $g(C_{i+1}) \leq
(d/d_{i+1} - 1)^2 $ as desired.    This procedure stops after
a finite number of  steps and  this completes the proof of  the
claim.
\enddemo

\qquad

\quad

\noindent{\bf \S 3 $d$-gonality of the  modular curve $X_0(N)$}

\vskip 0.5cm

\noindent In this section, we show our main theorem on
$d$-gonality of  the modular curves $X_0(N)_{\C}$ of level $N$.

The main idea of the proof obviously goes back to a  beautiful
argument of Ogg in [O],  which determines the complete list of
hyperelliptic modular curves.

For a positive integer $N$, let $\Gamma_0(N)$, $Y_0(N)_{\C}$  and
$X_0(N)_{\C}$ be as in Introduction.
  Classically, it is   well-known that $Y_0(N)_{\C}$  and
$X_0(N)_{\C}$ admit structures of algebraic  curves
$Y_0(N)_{\Q}$  and $X_0(N)_{\Q}$  defined over $\Q$, and
moreover by Igusa ([I]),
 there exist smooth models  of
$Y_0(N)$ and $X_0(N)$ over $\Z[1/N]$.   The smooth model
$Y_{0}(N)/\Z[1/N]$ can be considered as the coarse moduli space
classifying a pair $(E, C)$ of an elliptic curve $E$ with a cyclic
subgroup of order $N$ and  $X_{0}(N)/\Z[1/N]$  is its natural
compactification.  (See also [D-R] and [Ka-Ma] for more modern
treatment of  these facts.)
 Therefore for a prime $p$ with $p \nmid N$, $X_{0}(N)_{\Q}$  has
a good reduction at $p$.

Let $N$ be a positive integer and $p$ a prime such that $ p \nmid
N$ and $X_0(N)_{\F_p}$ the smooth projective curve over $\F_p$
obtained by a reduction of $X_{0}(N)_{\Q}$.  For any extension
field $k$ of
$\F_p$, a $k$-rational point of $X_0(N)_{k}$ corresponds to  an
isomorphism class of a pair $(E, C)$ of an elliptic curve with  a
cyclic subgroup of order $N$ defined over $k$.

The following argument is essentially due to Ogg ([O, Theorem
3]),  which gives  the lower bounds of number of
$\F_{p^2}$-rational points of
$X_0(N)_{\F_{p^2}}$.  (See also  [Ha-Si, Lemma 6].) If there
exists a supersingular elliptic curve $E$ defined over $\F_p$,
its Frobenius  endomorphism $\pi_p$ satisfies $\pi_p^2 = -p$.
Therefore any cyclic  subgroup $C \subset E$ of order $N$ is
defined over $\F_{p^2}$, which yields  a point $(E, C) \in
X_0(N)_{\F_{p^2}}(\F_{p^2})$.

For any prime $p$, set
$$  s(p) = \sum_{E/\F_p, \text{supersingular}}
\frac{1}{|\Aut(E)|}.
\tag 3.1
$$

Moreover  let $\nu(N)$ be the number of distinct prime factors of
$N$,  and set
$$
\phi(N) = [SL_{2}(\Z)/\pm 1:\Gamma_{0}(N)] = N \prod_{p|N} (1 +
1/p).
\tag 3.2
$$  From the argument above,  we can  obtain the following
lemma.   For detailed  proofs, see [O, Theorem 3] and  [Ha-Si,
Lemma 6].

\proclaim{Lemma 3.1} Let $p \nmid N$ be a prime.  Then we have
$$
\# X_{0}(N)_{\F_{p^2}} (\F_{p^2})  \geq  2^{\nu(N)} + 2s(p)
\phi(N).
\tag 3.3
$$  Moreover for $p = 2, 3$,  $s(2) =  1/24, s(3) = 1/12$, hence
we have
$$
  \# X_{0}(N)_{\F_{4}} (\F_{4}) \geq 2^{\nu(N)} + \frac{1}{12}
\phi(N), \quad \text{if}
\quad
  2 \nmid N, \tag 3.4
$$  and
$$
  \# X_{0}(N)_{\F_{9}} (\F_{9}) \geq 2^{\nu(N)} +
\frac{1}{6}\phi(N), \quad \text{if}
\quad
  3 \nmid N. \tag 3.5
$$
\endproclaim

\vskip 1cm

\proclaim{Theorem 3.2}  Let $X_0(N)$ be the modular curve of
level $N$.  If $X_0(N)_{\C}$  is d-gonal, that is, if it admits a
finite map
$f:X_0(N)_{\C} \lra \BP^1_{\C}$ of degree $d$,  then we have the
following.

1) If $N$ is odd, then we have
$$
\frac{1}{12}\phi(N) + 2^{\nu(N)}   \leq   \left\{ \matrix  5d  &
\text{if $d = 1, 2$}
 \\ 4(d -1)^2 + 5  & \text{if $d  \geq 3$}
\endmatrix \right. . \tag 3.6
$$

2) If $3 \nmid N$, then we have
$$
\frac{1}{6}\phi(N) + 2^{\nu(N)} \leq   \left\{ \matrix 10 d   &
\text{if $d = 1, 2$}
 \\  6(d -1)^2 + 10  & \text{if $d  \geq 3$}
\endmatrix \right. .\tag 3.7
$$
\endproclaim

\demo{Proof} From Tower Theorem 2.1,   we obtain a smooth
projective curve $C_{\Q}$ defined over $\Q$ and a finite morphism
$f:X_{0}(N)_{\Q} \lra C_{\Q}$ defined over
$\Q$ of degree $d'$ such that $ 1 \leq d' \leq d$  and  $g(C) \leq
(d/d' -1 )^2$.    First assume  that $g(C) \geq 1$.    For a
prime $p\nmid N$, the curve
$X_{0}(N)_{\Q}$ has  good reduction at $p$,  hence by Good
Reduction Lemma 5.1 (see \S 5),
$C_{\Q}$ has also good reduction at $p$ and we obtain a finite
morphism
$$ f \times \F_p: X_{0}(N)_{\F_p} \lra C_{\F_p}
$$ defined over $\F_p$.  Let $\F_q$ be the finite extension of
$\F_p$ with $q = p^2$.  Then since $g(C) \leq (d/d' -1 )^2$,  by
Weil's theorem of the analogue of the Riemann hypothesis, we can
bound the number of
$\F_q$-rational points of $C_{\F_q}$  as
$$
\# C_{\F_q}(\F_q) \leq 1 + 2 g(C) \sqrt{q} + q  = 1 + 2 p g(C) +
p^2.
\tag 3.8
$$ Since  $g(C) \leq (d/d' - 1)^2$,  from (3.8) above,  we obtain
$$
\# (X_{0}(N)_{\F_q} (\F_q)) \leq d' \cdot \# C_{\F_q}(\F_q)
\leq d' \cdot ( 1 + 2(d/d' -1)^2 p + p^2).
$$ Now fixing $d$,  set $H(p, d') =   d' \cdot ( 1 + 2(d/d' -1)^2
p + p^2)$.   For $1 \leq d' \leq d$, we can easily see that
$$ H(p, d') \leq \max \{ H(p, 1), H(p, d) \} = \max
 \{ p^2 + 1 + 2p \cdot (d -1)^2 ,  (p^2 + 1)d \} \tag 3.9
$$ From this and  Lemma 3.1,  putting  $p = 2, 3$, we obtain the
assertions (1) and (2). If $g(C) =0$, then we have a finite
morphism
  $X_{0}(N)_{\F_p} \lra C''_{\F_p}$ defined over
$\F_p$ but of  degree $1 \leq d" \leq d$ where
$C''_{\F_p}$ is a rational curve defined over
$\F_p$.  In this case, the bound becomes better than the former
case, which completes the  proof of theorem.

\enddemo

By using obvious inequalities $\phi(N) \geq N + 1$  and $\nu(N)
\geq 1$, we obtain the following corollary.

\proclaim{Corollary 3.3}  Let $X_0(N)$ be the modular curve of
level $N$.  If $X_0(N)_{\C}$  is d-gonal, that is, if it admits a
finite map
$f:X_0(N)_{\C} \lra \BP^1_{\C}$ of degree $d$,  then we have the
following.

1) If $N$ is odd, then we have
$$ N  \leq   \left\{ \matrix  60 d - 25   & \text{if $d = 1, 2$}
 \\  48(d - 1)^2 + 35  & \text{if $d  \geq 3$}
\endmatrix \right. .\tag 3.10
$$

2) If $3 \nmid N$, then we have
$$ N  \leq   \left\{ \matrix 60d - 11   & \text{if $d = 1, 2$}
 \\  36(d - 1)^2 + 47  & \text{if $d  \geq 3$}
\endmatrix \right. .\tag 3.11
$$
\endproclaim

\proclaim{Remark 3.4}
{\rm Since we know the genus formula of
$X_{0}(N)_{\C}$, we know all the cases with $ g(N) := g(X_{0}(N)_{\C})
\leq 1$.  If $g(N) = 0$, we have $N = 1, \cdots, 10, 12, 13, 16, 18, 25$, and
if $g(N) = 1$, we have $N = 11, 14, 15, 17, 19, \cdots, 21, 24, 27, 32, 36,
49$.
Ogg [O] showed that $X_{0}(N)_{\C}$ is a hyperelliptic
curve of $g(N) \geq 2$ if and
only if
$N = 22, 23, 26, 28, 29, 30, 31, 33, 35, 37,  39, 41, 46, 47, 50, 59,
71$,  (i.e. 19 values).  This implies that if $X_0(N)_{\C}$
is 2-gonal then $N \leq 71$.  This is obviously better than our
bounds in Corollaries 3.3 and 3.5, because we use only  rough
estimates of $\phi(N)$ and $\nu(N)$.}
\endproclaim

\vskip  0.5cm

For a general positive integer $N$, write $N = 2^l \cdot M$ such that
$M $ is odd.  Then we have natural finite morphisms

$$
\varphi_1:X_{0}(N)_{\C} \lra X_{0}(2^l)_{\C},
 \quad \varphi_2:X_{0}(N)_{\C} \lra X_{0}(M)_{\C}.
$$
 By Lemma 1.3, if $X_{0}(N)_{\C}$ is $d$-gonal, then  both of
$X_{0}(2^l)_{\C}$ and $X_{0}(M)_{\C}$ are
$d$-gonal.  Hence we have the following corollary  which gives a
bounds for $N$ by  a polynomial  in $d$  of degree $\leq 4$.

\proclaim{Corollary 3.5}  Let $X_0(N)$ be the modular
curve of level $N$.  If $X_0(N)_{\C}$  is d-gonal,
that is, if it admits a finite map
$f:X_0(N)_{\C} \lra \BP^1_{\C}$ of degree $d$,  then we have the
following.
$$
N  \leq   \left\{ \matrix  (60 d -25)(60 d -11)   & \text{if $d = 1, 2$}
 \\  (48(d - 1)^2 + 35)( 36(d - 1)^2 + 47 ) & \text{if $d  \geq 3$}
\endmatrix \right. . \tag 3.12
$$
\endproclaim

\proclaim{Remark 3.6}
{\rm Let $k = \F_p$ and fix an algebraic closure $\k$ of $k$.
For $N$ such that $ p \nmid N$, assume that the smooth
projective curve
$X_0(N)_{\k}$ is d-gonal.  From Tower Theorem 2.1,
we obtain a smooth projective curve $C'$ defined over $k$ and
a finite morphism $f:X_{0}(N)_{k} \lra C'$  of degree
$d', 1 \leq d' \leq d$ satisfying that
$ g(C') \leq (d/d' -1)^2$.  From the same arguments
as in Theorem 1.2, we obtain an inequality
$$
\# (X_0(N)_{\F_{p^2}}(\F_{p^2}) \leq \max \{p^2 + 1 + 2p \cdot (d -1)^2,
(p^2 + 1) \cdot d \}.
$$
Then from Lemma 3.1, (3.3) , we obtain an inequality
$$
2^{\nu(N)} + 2 s(p) \phi(N) \leq \max \{p^2 + 1 + 2p \cdot (d -1)^2,
(p^2 + 1) \cdot d \}.
$$
Hence if $s(p) > 0$, then we obtain a bound for $N$ by a constant
only depending on $d$ and $p$.}
\endproclaim

\vskip 1cm

\noindent {\bf \S 4 Bounding  orders of torsions    of
Mordell-Weil group.}

\quad

   Let $C$ be a  smooth projective curve defined over an
algebraically closed field $k$ with the rational function  field
$K = k(C)$, and let $E$ be an elliptic curve defined over
$K$.  Then we obtain
 a relatively minimal elliptic surface $\pi:{\Cal E}
\lra C$  associated  to $E/K$.  We call $E/K$ is ``constant'', if
its $K/k$-trace is  non-trivial.   If $E/K$ is constant, then the
associated family
$\pi:{\Cal E}
\lra C$ is birational equivalent to $E_{0} \times C$ with  an
elliptic curve $E_{0}$ over $k$.
 It is known that    the Mordell-Weil group $E(K)$ of $E$ is
finitely  generated,  if $E/K$ is not constant (cf. [La]).

\proclaim{Lemma 4.1}  Let $C$ and $E/K$ be as above, and assume
that the characteristic of the base field $k$ is zero and
$E/K$ is not constant.   Then if the Mordell-Weil group
$E(K)$ has a cyclic subgroup of order $N > 1$.  Then one obtains
a surjective morphism
$h:C \lra X_0(N)_{\C}$.
\endproclaim

\demo{Proof}
 Let $\pi^0: {\Cal E}^0 \lra C^0$ denote  the morphism   obtained
by  restricting $\pi$ to the maximal  open set ${\Cal E}^0 $ on
which  $\pi$ is smooth.   Now assume that  the Mordell-Weil group
$E(K)$ has a cyclic subgroup of order
$N$, it defines a cyclic group of order
$N$  on each fiber $\Cal E_t$  for
$t \in C^0$.    Then since $Y_0(N)_{\C}$ is the moduli space of
pairs $(E, D)$ (see [D-R] and [Kat-Ma]),   we have a natural
non-constant morphism
$h^0:C^0 \lra Y_0(N)$,  which extends to a finite  morphism
$h: C \lra X_0(N)$.
\enddemo

{}From Lemma 4.1 together with Lemma 1.3, we have the following

\proclaim{Proposition 4.2}   Let $C$ be a  smooth projective curve
defined over $k$ with the rational  function field $K$
 and assume that $C$ is
$d$-gonal.  If  there exists a non-constant  elliptic curve $E$
over $K$ whose Mordell-Weil  group
$E(K)$ has a cyclic group of order $N$, then  the modular curve
$X_0(N)$ is also $d$-gonal.
\endproclaim

Together with Corollaries 3.3 and 3.5,  this proposition implies
the following theorem, which may be considered as  an analogue of
strong uniformly boundedness conjecture  in the function field
case.

\proclaim{Theorem 4.3}  Let  $k$ be an algebraically closed field
of characteristic zero and $C$ be a smooth projective curve
defined over $k$,  and let $K = k(C)$ be the function field of
$C$. Then if there exists a non-constant elliptic curve $E$ defined over $K$
such that its Mordell-Weil group $E(K)$ admits a torsion element of
order $N$.  Then we have a polynomial function $B(d)$ in $d$
such that
$$
  N \leq B(d). \tag 4.1
$$
For example, we take B(d) as the right hand side of (3.12).

\endproclaim

\vskip 1cm

\noindent {\bf \S 5  A lemma of good reduction of morphisms}

\vskip 0.5cm

In this section, we shall use the following lemma, which we have used in
the proof of Theorem 3.2.

\proclaim{Lemma 5.1} Let $C_1$ and $C_2$ be projective  smooth
curves defined over $\Q$  both of which are  geometrically
irreducible,  and let $f:C_1
\lra C_2$ be  a dominant morphism of degree $d$ which  is also
defined over $\Q$.  Assume that $C_1$ has good reduction at a
prime integer $p  > 0$.  Then we have the following.

1)  Assume that $g(C_1) \geq 1$.   Then $C_2$ has  a good
reduction  at $p$ and $f$ induces a finite morphism of degree
$d$
$$ f_s: C_{1, s}  \lra C_{2,s},
$$  defined over $\F_p$  where $C_{i, s}$ $ i =1, 2$ denote  the
smooth projective  curves defined over
$\F_p$ obtained from $C_i$ respectively.

2) If $g(C_2) = 0 $, then we obtain a dominant morphism
$$ f^{'}_{s}:C_{1,s} \lra  C'
$$ of degree $d'  \leq d$ defined over $\F_p$ such that
$C'$ is a smooth rational curve defined over $\F_p$.
\endproclaim

\demo{Proof} Let us set $S = Spec(\Z_p)$ where $\Z_p$ denotes the
ring of  p-adic integers, and denote by
$\eta$ and $s$ the generic and  the closed points of
$S$ respectively.  Moreover, we set
$\hat{S}  =Spec(\Z_p^{sh})$, where $\Z_p^{sh}$ is the strict
henselization of
$\Z_p$ whose residue field $k = \Z_p^{sh}/(p)$ is a fixed algebraic
closure $\ovl{\F_p}$ of
$\F_p = \Z_p/(p)$.   For each curve $C_i$ over $\Q$, we obtain a
regular $\Z_p$-model of
$C_i$
$$
\pi_i:\calC_i \lra S = Spec(\Z_p)
$$ such that:

i) $\calC_i$ is a regular scheme,

ii) $\pi_i$ is a proper flat morphism and

iii)  $\calC_{\eta} \simeq C_i \times_{\Q} \Q_p$.

Such models can be obtained by the blowing up of any projective
model of $C_i$ over $S$ by blowing up (cf. [Lip].) Moreover since
$C_1$ has good reduction, we can assume that
$\pi_1:\calC_1 \lra S$ is smooth.

Now let us assume that $g(C_2) \geq 1$,  Then, by Lichtenbaum
[Lic] and Shafarevich [Sha], we may  assume that
$\pi_2:\calC_2 \lra S$ is the minimal model.  Let
$J_i, i =1, 2$ denote the Jacobian varieties of $C_i$
respectively.   Then since
$J_1$ has good reduction at $p$,  it is easy to see that $J_2$ also
has good reduction at $p$.  (For a proof, one may use Serre-Tate
criterion [S-T, Cor. 2], or just use the N\'{e}ron models of $J_1$ and
$J_2$.)
   Since there exists a dominant morphism $C_1 \lra C_2$ defined
over $\Q$, we obtain  a dominant
$S$-rational map
$$
\varphi:\calC_1 \cdots \rightarrow \calC_2.
$$
Now  considering the base extension $\hat{S} \lra S$,  we obtain
the induced
$\hat{S}$-rational map
$$
\matrix
   &  &  \varphi  &  &   \\
\ovl{\calC_1}  & &\cdots \rightarrow & &\ovl{\calC_2}
\\
  \pi_1& \searrow &   & \swarrow & \pi_2\\
   &  &  \hat{S} &    &
\endmatrix
$$
Note that $\ovl{\calC_2}$ is also the minimal model over
$\hat{S}$.   Since $J_2$ has  good reduction at $p$,  it implies
that
$\ovl{\calC_2}$ has at least semistable reduction at
$p$ by [D-M, Theorem 2.4]. Consider the irreducible decomposition of $
\ovl{\calC_{2,s}}$:
$$
\ovl{\calC_{2,s}} = \sum_{i=1}^l T_i.
$$
The dual graph of $\ovl{\calC_{2,s}}$  has no cycle because  the
N\'eron model has no torus part, hence each irreducible component
$T_i$ is smooth and because $J_2$ has good reduction at $p$,  one
has $\sum_{i=1}^l g(T_i) = g(C_2) \geq 1$. Now we  claim that:

\medskip {\bf Claim:} {\it  $\ovl{\calC_{2,s}}$ has only one
irrational component, say,
$T_1$. }
\medskip

If we admit this claim, we can show that
$\ovl{\calC_{2,s}} = T_1$, which  also implies that
$\ovl{\calC_{2,s}}$ is smooth.
Let $l$ denote the the number of
irreducible components of
$\ovl{\calC_{2,s}}$ and
assume that $l > 1$.  Then since the dual graph
of
$\ovl{\calC_{2,s}}$ has no cycle, there  exists a component
$T_i$, $i \geq 2$ such that
$$
(\ovl{\calC_{2,s}} - T_i) \cdot T_i = 1.
$$
 Since $\ovl{\calC_{2,s}} \cdot T_i = 0$, we have
$(T_i)^2 = -1$.   However since all  components $T_j$,
$ j \geq 2$  are smooth rational curves, this implies that $T_i$ is
an exceptional  rational curve of the first kind, which contradicts
to the minimality of  the model $\ovl{\calC_2}$.

Now  it is easy to see that
$\pi_2:\ovl{\calC_2} \lra \hat{S}$ is smooth, and since $\hat{S}
\rightarrow  S$ is faithfully flat, we conclude that
$\pi_2$ is smooth.

Now we  prove the claim. After a sequence of quadratic
transformations with centers in the  special fiber of
$\pi_1$, we obtain  a birational morphism
$\tau:\ovl{\Gamma} \lra
\ovl{\calC_1}$
such that $\hat{S}$-rational map
$\varphi$ induces a  dominant
$\hat{S}$-morphism $\phi:\ovl{\Gamma} \lra
\ovl{\calC_2}$. Note that since $\pi_1$ and
$\pi_2$ are proper $\phi$ must  be surjective. Consider the
following commutative diagram (cf. [Sha]):
$$
\matrix
    & & \ovl{\Gamma} &  &  \\
\tau &  \swarrow &  & \searrow \phi &\\
   &  &  \varphi  &  &   \\
\ovl{\calC_1}  & &\cdots \rightarrow & &\ovl{\calC_2}
\\
  \pi_1& \searrow &   & \swarrow & \pi_2\\
   &  &  \hat{S} &    &
\endmatrix
$$
Set $F = \ovl{\calC_{1, s}}$, the special fiber of
$\pi_1$  and let $F'$ be the proper transform of $F$  by
$\tau$,
$\ovl{\Gamma}_s$  the special fiber of $\tau \circ
\pi_1$, and  write
$$
\ovl{\Gamma}_s = F' + \sum_{i} m_iE_i.
$$
Here all of  $\{E_i\}$ are exceptional rational curves obtained
by  the blowing up
$\tau$.   Now consider the morphism
$\phi_s:\ovl{\Gamma}_s \lra
\ovl{\calC_{2,s}}$. Since
$F'$ is the only one irrational component of
$\ovl{\Gamma}_s$,
$\ovl{\calC_{2,s}}$  has at most one irrational component, which
proves the claim.

Next we  prove  assertion 2) of the lemma.  Assume that
 $g(C_2) = 0$ and let
$$
f:C_1 \lra C_2
$$
be the dominant morphism of degree $d$.

 Then we obtain a dominant
$\Z_p$-rational map
$f: \calC_1 \cdots \rightarrow \calC_2$.  As in the former
case,  by Shafarevich [Sha], we obtain the commutative diagram
$$\matrix
    & & \Gamma &  &  \\
\tau &  \swarrow &  & \searrow \phi &\\
   &  &  f  &  &   \\
\calC_1  & &\cdots \rightarrow & & \calC_2 \\
  \pi_1& \searrow &   & \swarrow & \pi_2\\
   &  &  S. &    &
\endmatrix
$$
Here $\tau:\Gamma \rightarrow \calC_1$ is obtained by  a sequence
of quadratic  transformations with centers  in the  special fiber of
$\pi_1$.  Let $F'$ be the proper  transform of
$F = \calC_{1, s}$.  If the restriction of $\tau$ to
$F'\simeq \calC_{1, s}$ is not constant, the image
$\phi(F')$ must be  a rational curve
$\calC_{2, s}$. The degree of the obtained map
$F' \lra
\calC_{2, s}$ is equal to  or  less than $d $.   If the restriction
of $\phi$ to
$F'$ is constant, we will blow up $\calC_2$ at  the closed
point  $x = \tau(F')$ and its infinitely near points  in order to
make the morphism $\tau$ flat.  As mentioned in remark before the
proof of [B-L-R, Prop. 6, 3.5] (see also [R-G]),  there exists a blowing up $V
\lra
\calC_{2, s}$ such that  the induced morphism
$\phi':\Gamma' \lra V$ is flat.  (Here
$\Gamma'$ is the schematic closure of  $\Gamma_{\eta}$ in
$\Gamma \times_{\calC_2} V$.)
  Since $\phi'$ is flat, it maps the  proper transform of $F'$
$(\simeq \calC_{1,s})$
 by $\Gamma' \rightarrow \Gamma$ onto a curve $C'$ which is a
rational curve over $\F_p$ arising as an
 exceptional curve of the blowing up.  Therefore, we obtain  a
surjective morphism
$\calC_{1, s} \lra C'$, whose degree $d'$ is less than  or equal to
$d$, the degree of $f$ at the generic fiber.
\enddemo

\proclaim{Remark 5.2}{\rm Without changing the above proof, we can extend
the statement of Lemma 5.1
to the case of curves over a  discrete valued field $K$ with the
integer ring $R$ under the assumption that the residue field $R/m$ is
perfect.}
\endproclaim

\enddocument